\documentclass[12pt,final]{article}
\usepackage{amsmath,amsthm,amsfonts,amssymb,graphicx,enumerate,psfrag}
\usepackage{tikz,pgfplots}
\usepackage{mathrsfs}
\usepackage{fullpage}

\usepackage[colorlinks,citecolor=blue,urlcolor=blue]{hyperref}
\definecolor{myblue}{RGB}{51,51,178}
\definecolor{myred}{RGB}{189,26,26}
\definecolor{mygreen}{RGB}{0,128,0}
\usepackage[utf8]{inputenc} 
\newtheorem{lemma}{Lemma}

\newtheorem{theorem}{Theorem}

\newtheorem{proposition}{Proposition}
\newtheorem{corollary}{Corollary}

\newcommand{\dN}{\mathbb {Z}_+}
\newcommand{\dZ}{\mathbb {Z}}

\newcommand{\cU}{\mathcal {U}}

\newcommand{\cP}{\mathcal {P}}
\newcommand{\bI}{\mathbb{I}}
\newcommand{\cL}{\mathcal {L}}

\newcommand{\dR}{\mathbb {R}}
\newcommand{\dd}{\mathrm{d}}

\newcommand{\EE}{{\mathbb{E}}}
\newcommand{\PP}{{\mathbb{P}}}

\newcommand{\tmix}{\mathrm{t}_{\textsc{mix}}}
\newcommand{\wmix}{\mathrm{w}_{\textsc{mix}}}

\newcommand{\dtv}{\mathrm{d}_{\textsc{tv}}}

\newcommand{\Ent}{{\mathrm{Ent}}}

\newcommand{\Var}{{\mathrm{Var}}}

\newcommand{\gap}{{\mathrm{gap}}}
\newcommand{\wnu}{{\widehat{\mu}}}
\newcommand{\Varent}{{\mathrm{Varent}}}
\newcommand{\supp}{\mathrm{supp}}
\title{Cutoff on abelian groups}
\author{Guy Blachar, Jonathan Hermon and Justin Salez}
\begin{document}
\maketitle
\begin{abstract}
We prove the long-predicted equivalence between the cutoff phenomenon and the product condition along any sequence of abelian random walks whose inverse transition probabilities grow at most logarithmically with the size of the group. Our result does not require reversibility, and comes with a sharp quantitative estimate on the mixing window. Moreover, we show that our  sparsity restriction is optimal by constructing abelian Cayley graphs which satisfy the product condition but fail to exhibit cutoff, and whose degrees are arbitrarily close to being logarithmic in the volume.
\end{abstract}
\tableofcontents
\section{Introduction}
\subsection{Main estimate}
Throughout this paper, $G$ is a finite abelian group equipped with a probability measure $\mu\in\cP(G)$, and we consider a $G$-valued continuous-time Markov chain $X=(X_t)_{t\ge 0}$ with initial condition $X_0=0$ and generator $\cL$ acting on $L^2(G)$ as follows:
\begin{eqnarray}\label{def:L}
\cL f(g) & := & \sum_{s\in G} \mu(s)\left[f(g+s)-f(g)\right].
\end{eqnarray}
We call $X$ an \emph{abelian random walk} with jump distribution $\mu$. 
An important special case is that where $\mu$ is uniform on a set $S\subseteq G$:  the process $X$  is then exactly a simple random walk on the \emph{Cayley graph} $\Gamma=\mathrm{Cay}(G,S)$, i.e. the directed graph with vertices and edges
\begin{eqnarray}
V_\Gamma\ := \ G, & \textrm{ and } & E_\Gamma:=\{(g,g+s)\colon g\in G,s\in S\}.
\end{eqnarray}
Coming back to the general case, we may, upon restricting the group if needed, assume  that  $G$  is generated by the support of $\mu$. This guarantees  that the law of $X_t$ approaches the uniform law $\cU$ on $G$ as $t\to\infty$. A standard way to quantify this convergence is through the \emph{total-variation distance to equilibrium}, defined at any time $t\ge 0$ as follows:
\begin{eqnarray}
\dd(t) & := & \max_{A\subseteq G}\left|\PP\left(X_t\in A\right)-\cU(A)\right|. 
\end{eqnarray}
This function  decreases continuously from $\dd(0)=1-1/|G|$ to $\dd(\infty)=0$, and the time at which it reaches a given precision $\varepsilon\in(0,1)$ is known as the \emph{mixing time}:
\begin{eqnarray}
\tmix(\varepsilon) & := & \inf\{t\ge 0\colon \dd(t)\le \varepsilon\}.
\end{eqnarray}
The default value for $\varepsilon$ is $1/(2e)$, in which case we simply write $\tmix$. This quantity describes \emph{how long} the process needs to be run in order to approach equilibrium. Motivated by the celebrated \emph{cutoff phenomenon}, we here seek to understand \emph{how abrupt} the transition from out-of-equilibrium to equilibrium is. This is captured by the difference
 \begin{eqnarray}
\wmix(\varepsilon) & := & \tmix(\varepsilon)-\tmix(1-\varepsilon),
\end{eqnarray}
known as the \emph{width} of the \emph{mixing window}. Our main contribution is a sharp estimate on this subtle quantity, in terms of two basic statistics of the random walk. The first one controls the \emph{sparsity} of the jump distribution $\mu$, through the smallest non-zero jump rate:
\begin{eqnarray}
\mu_\star & := & \min\left\{\mu(s)\colon s\in\supp(\mu)\right\}.
\end{eqnarray}
The second is the \emph{spectral gap} of $\cL$,  i.e. the difference between its two largest eigenvalues, once projected onto the real axis. Since $\cL$ commutes with its transpose $\cL^\top$, this is also the spectral gap of the symmetric generator $(\cL+\cL^\top)/2$, and we have the variational characterization
\begin{eqnarray}
\gap & = & \inf\left\{{\langle f,-\cL f\rangle}\colon f\in L^2(G),\langle f,1\rangle =0,\langle f,f\rangle=1\right\}.
\end{eqnarray}
With this notation at hand, our main result is the following simple  estimate on the width of the mixing window, valid for any abelian random walk. 
\begin{theorem}[Width of the mixing window]\label{th:wmix}For any $\varepsilon\in(0,1/2)$, we have
\begin{eqnarray}
\label{wmix}
\wmix(\varepsilon) & \le & C_\varepsilon\sqrt{\frac{\tmix}{\gap}\max\left\{1,\log\left(\frac{1}{\mu_{\star}\tmix}\right)\right\}},
\end{eqnarray}
where $C_\varepsilon<\infty$ is a constant that depends  only on $\varepsilon$.
\end{theorem}
This  estimate will be shown to be sharp, except for the constant $C_\varepsilon$ (see Theorem \ref{th:sharp} below). It  becomes particularly interesting when the size of the group $G$ gets large, because it implies that certain abelian random walks exhibit a remarkably abrupt transition from out-of-equilibrium to equilibrium known as a \emph{cutoff phenomenon}. 
\subsection{Implications for cutoff}
As is often the case in concrete models, we now assume that the pair $(G,\mu)$  depends on a size parameter $n\ge 1$, which is kept implicit   to lighten notation, and we consider the large-volume limit where $n\to\infty$. In this regime, it is sometimes the case that
\begin{eqnarray}
\label{def:cutoff}
\forall \varepsilon\in(0,1/2),\quad \tmix & \gg & \wmix(\varepsilon).
\end{eqnarray}
Here and throughout the paper, we use the standard asymptotic notation  $u\gg v$, $u\sim v$,  $u=\Omega(v)$, $u=O(v)$ and $u=\Theta(v)$ to mean, respectively, that the ratio $u/v$ tends to $\infty$,  to $1$, remains bounded away from $0$, from $\infty$, and from both as $n\to\infty$.  In words, (\ref{def:cutoff}) means that instead of decaying gradually from $1$ to $0$ over time, the distance to equilibrium $t\mapsto\dd(t)$ approaches a step function as $n\to\infty$. Such an abrupt transition from out-of-equilibrium to equilibrium is known as a \emph{cutoff}. We refer the reader to the seminal works \cite{DiaconisShahshahani1981,Aldous1983,AldousDiaconis1986,Diaconis1996} and the recent lecture notes \cite{salez2025lecture} for an introduction to this fascinating phenomenon. Despite the accumulation of several examples, an exact characterization is still missing, even in the simple setting of abelian random  walks. As a corollary of Theorem \ref{th:wmix}, we obtain the following  simple sufficient condition for a sequence of abelian random walks to exhibit cutoff. 

\begin{corollary}[Cutoff]\label{co:cutoff}A sequence of abelian random walks exhibits cutoff as soon as
\begin{eqnarray}
\gap\, \tmix & \gg & \max\left\{1,\log\left(\frac{1}{\mu_{\star}\tmix}\right)\right\}.
\end{eqnarray}
\end{corollary}
To the best of our knowledge, this criterion covers all instances of cutoff on abelian groups in the literature, including those where the jump distribution is random  \cite{Hough2017,HermonOleskerTaylor2026}. It also provides  new interesting examples, such as the abelian sandpile Markov chain introduced in \cite{kaiser2024abelian}, which has $\tmix=\Theta(n\log n)$, $\gap=\Theta(1/n)$, $\mu_\star=\Theta(1/n)$, and hence exhibits cutoff. Perhaps more importantly, Corollary \ref{co:cutoff} has the following conceptual interest. Recall that any sequence of Markov chains exhibiting cutoff must  satisfy the celebrated \emph{product condition}
\begin{eqnarray}
\label{assume:product}
\gap\, \tmix  & \gg & 1,
\end{eqnarray}
but that this is generically not enough; see, e.g., \cite[Section 1.5]{salez2025lecture}.  For abelian random walks, Corollary \ref{co:cutoff}  shows that this classical \emph{necessary} condition is actually  \emph{sufficient}, as long as
\begin{eqnarray}
\label{assume:good}
\mu_\star\,\tmix & = & \Omega(1).
\end{eqnarray}
In words, the condition (\ref{assume:good}) simply means that each generator $s\in\supp(\mu)$ has a non-vanishing chance of being actually used during mixing. As we will see, this natural requirement holds in particular throughout the \emph{logarithmically sparse} regime, defined as follows:
\begin{eqnarray}
	\label{assume:sparse}
	\frac{1}{\mu_\star} & = & O\left(\log|G|\right).
\end{eqnarray}
This leads to the following result, which  partially resolves \cite[Conjecture 1]{HermonHuangSalez2026} and adds the class of logarithmically sparse abelian random walks to the slowly growing list of Markov processes for which cutoff reduces to the product condition \cite{DingLubetzkyPeres2010,BasuHermonPeres2017,Salez2023Exclusion,Salez2026Diffusions,HermonHuangSalez2026}. 
\begin{corollary}[Logarithmic sparsity]\label{co:sparse}The cutoff phenomenon (\ref{def:cutoff}) is equivalent to the product condition (\ref{assume:product}) along any sequence of abelian random walks satisfying (\ref{assume:sparse}).
\end{corollary}
In particular, this completely characterizes cutoff for simple random walk on abelian Cayley graphs whose degrees are at most logarithmic in the volume. 
Finally, our next result shows that this restriction is optimal, and that our main window estimate  is sharp.
\begin{theorem}[Sharpness]\label{th:sharp}Let $(a_n)_{n\ge 1}$ be positive numbers satisfying $1\ll \log a_n\ll \log n$. Then, there  is a sequence of (symmetric) abelian Cayley graphs $\{\mathrm{Cay}(G_n,S_n)\}_{n\ge 1}$   such that
\begin{eqnarray}
|S_n| & \sim & a_n\log |G_n|,
\end{eqnarray}
on which simple random walk satisfies the product condition but fails to exhibit cutoff. Moreover,  the two sides of the window estimate (\ref{wmix}) are of the same order for small enough $\varepsilon>0$.
\end{theorem}
\subsection{Related work}

Random walks on abelian groups have long provided a rich class of models with enough structure and symmetry to serve as an ideal testbed for the exploration of mixing times and the cutoff phenomenon. A pioneering work of Diaconis and Saloff-Coste \cite{DiaconisSaloffCoste1994} implies, in particular, that cutoff
cannot occur when the number of generators remains bounded, due to the failure of the product condition. On the positive side, cutoff has been established for various  families of abelian
Cayley graphs with growing degrees, and most notably \emph{random instances} obtained by choosing the generators uniformly at random. This body of work builds on a fertile research program proposed by Aldous and Diaconis in  \cite{AldousDiaconis1986}, and on the famous discovery, by Alon and Roichman, that random Cayley graphs exhibit optimal expansion  \cite{AlonRoichman1994}. See, in particular, the inspiring study \cite{Hough2017} by Hough  on cyclic groups, and the impressive series \cite{HermonOleskerTaylor2023,HermonOleskerTaylor2026} by Hermon and Olesker-Taylor on general abelian groups. All those instances are covered by our main result. 

A complementary line of research consists in seeking systematic criteria for the occurrence of 
cutoff phenomena, without relying on the detailed, model-specific  analysis of mixing times. In particular, the varentropy approach introduced by Salez in \cite{Salez2024JEMS} and  refined by Pedrotti and Salez in \cite{PedrottiSalez2025} shows that cutoff follows from the product-like condition
\begin{eqnarray}
	\label{prior}
\gap\,\tmix & \gg & \log \frac{1}{\mu_\star},  
\end{eqnarray}
along any sequence of Markov chains with non-negative curvature and symmetric support. (see, e.g., \cite[Theorem 5.3]{salez2025lecture}). When specialized to abelian random walks, this requirement is, however, much stronger than the one in our Corollary \ref{co:cutoff}. In particular, it is not even satisfied on the standard hypercube.  The present paper can be seen as an optimal refinement of the varentropy approach, when tailored to the special case of abelian random walks. 

Finally,  the relevance of abelian random walks extends beyond abelian
groups themselves. Indeed, the mixing behavior of several non-commutative random walks has been found to be  governed by that of their abelianization. Instances of this general \emph{dimensional reduction}  appear notably in the analysis of random walks on upper-triangular matrices by Arias-Castro, Diaconis and Stanley \cite{AriasCastroDiaconisStanley2004}, Nestoridi \cite{Nestoridi2019}, and Nestoridi and Sly \cite{NestoridiSly2023}, as well as in the study of random Cayley graphs by Hermon and Olesker-Taylor \cite{HermonOleskerTaylor2019}. More recently, Huang \cite{Huang2026Conjugacy}  and Hermon, Huang and Salez 
\cite{HermonHuangSalez2026} showed that a large class of random walks on nilpotent groups exhibit cutoff if and only if their abelianizations do.  This clearly motivates a thorough study of cutoff on abelian groups, which is the aim of the present paper.

\section{Proof of the main estimate}
Our proof of Theorem \ref{th:wmix} relies on the general varentropy approach developed in the series of works \cite{Salez2024JEMS,HermonHuangPedrottiSalez2026,Salez2026Diffusions,PedrottiSalez2025,salez2025lecture}. Let us briefly recall it in the present context, where the state space is finite and the stationary distribution uniform.
\subsection{The varentropy approach}
 The Shannon \emph{entropy} of a   $G$-valued random variable $Z$ with law $\nu$ is classically defined as
\begin{eqnarray}
\Ent(Z) & := & \sum_{g\in G}\nu(g)\log \frac{1}{\nu(g)} \ = \ \EE\left[\log\frac{1}{\nu(Z)}\right].
\end{eqnarray}
To define the \emph{varentropy} of $Z$, we simply replace the above expectation with a variance:
\begin{eqnarray}
\Varent(Z) & := &  \Var\left[\log\frac{1}{\nu(Z)}\right].
\end{eqnarray}
Following the lecture notes \cite{salez2025lecture}, we say that the random walk $(X_t)_{t\ge 0}$ satisfies an \emph{information-theoretic differential inequality} (IDI) with rate $\Psi\colon (0, \infty) \to (0, \infty)$ if
\begin{eqnarray}
\label{def:IDI}
\forall t> 0,\quad \Varent(X_t) & \le & \Psi(t)\,\frac{\dd}{\dd t}\Ent(X_t). 
\end{eqnarray}
As shown in   \cite[Lemma 5.4]{salez2025lecture}, this  inequality systematically implies the window estimate
\begin{eqnarray}
\label{IDI:wmix}
\forall \varepsilon\in(0,1/2),\quad \varepsilon^3\wmix(\varepsilon) & \le & \frac{1}{\gap}+\sqrt{\frac{M_\varepsilon}{\gap}},
\end{eqnarray}
where $M_\varepsilon$ denotes the supremum of the rate function $\Psi$ over the interval $[\tmix(1-\varepsilon),\tmix(\varepsilon)]$. This approach was first implemented in  the context of non-negatively curved diffusions \cite{Salez2026Diffusions}, which were shown to satisfy the IDI (\ref{def:IDI}) with  $\Psi(t)=2t$. In view of (\ref{IDI:wmix}), this implies that their cutoff is characterized by the product condition. The theory was later extended to non-negatively curved Markov chains on finite state spaces  \cite{PedrottiSalez2025}, including abelian random walks. However,  the lack of a chain rule on discrete spaces incurs an extra logarithmic dependency in  $\frac{1}{\mu_\star}$ (see (\ref{prior}) above), and the latter always diverges in the cutoff regime \cite{MunchSalez2023}. Using a different approach, we here manage to replace the term $\frac{1}{\mu_\star}$ with the much smaller (and optimal) quantity $\frac{1}{\mu_\star \tmix}$, which crucially remains bounded throughout the regime (\ref{assume:sparse}). 
\begin{theorem}[Universal IDI]\label{th:IDI}Any abelian random walk satisfies the IDI (\ref{def:IDI}) with rate 
\begin{eqnarray}
\label{def:psi}
\Psi(t) & := & 550\, t \max\left\{1,\log \frac{1}{\mu_\star t}\right\}.
\end{eqnarray}
\end{theorem}
Before proving this result, let us quickly explain why it implies our main estimate.
\begin{proof}[Proof of Theorem \ref{th:wmix}] The function  $t\mapsto \dd(t)$ classically satisfies $\dd(t+s)\le 2\dd(t)\dd(s)$ for all $t,s\ge 0$, see, e.g., \cite[Lemma 1.1]{salez2025lecture}. Recalling that $\tmix=\tmix(1/(2e))$, this implies
\begin{eqnarray}
\tmix(\varepsilon) & \le & \left\lceil \log\frac{1}{\varepsilon}\right\rceil\,\tmix,
\end{eqnarray}
for any $\varepsilon\in(0,1)$. Since the function $\Psi$ defined at (\ref{def:psi}) is increasing, we immediately deduce that its supremum over the interval $[\tmix(1-\varepsilon),\tmix(\varepsilon)]$ satisfies
\begin{eqnarray*}
 M_\varepsilon & \le & 550 \left\lceil \log\frac{1}{\varepsilon}\right\rceil\tmix\max\left\{1,\log\frac{1}{\mu_\star \tmix}\right\}.
 \end{eqnarray*} 
On the other hand, the classical inequality 
$
\gap\,\tmix \ge 1
$ (see, e.g., \cite[Remark 1.1]{salez2025lecture}) ensures that $\frac 1 \gap\le \sqrt{\frac{\tmix}{\gap}}$. Thus, (\ref{IDI:wmix}) yields the desired conclusion, with $C_\varepsilon=\frac{1+\sqrt{550 \left\lceil \log\frac{1}{\varepsilon}\right\rceil}}{\varepsilon^3}$.
\end{proof}
In order to prove Theorem \ref{th:IDI}, we will exploit the  elementary observation that an abelian random walk can be realized as the  image, under a suitable projection, of independent Poisson processes (namely, those counting the jumps in each direction $s\in\supp(\mu)$). As a consequence, varentropy estimates along the random walk can be deduced from appropriate \emph{weighted Poincaré inequalities} for Poisson distributions. 
\subsection{Poincaré inequality for Poisson distributions}
Our starting point is a simple  lemma, which appears to be new and of independent interest.
\begin{lemma}[Poincaré inequality for attractive birth-death chains]\label{pr:BDC}Consider an irreducible birth-death  chain on an interval $\bI\subseteq\dZ$ with non-increasing birth rates $(q_+(x))_{x\in\bI}$, non-decreasing death rates $(q_-(x))_{x\in\bI}$ and stationary law $\pi$.  Then, the Poincaré inequality
\begin{eqnarray}
\label{BDCPoincaré}
\forall f\in L^2(\pi),\quad 
\Var_{\pi}\left(f\right) & \le &  C\,\sum_{x\in  \bI}\pi(x)q_+(x)\left[f(x+1)-f(x)\right]^2,
\end{eqnarray}
holds with $C = \inf_{o\in\bI} \frac{4}{\pi^2(o)q_\star(o)}$,  where $q_\star(o)$ denotes the smallest non-zero jump rate to $o$.
\end{lemma}
\begin{proof}
A classical use of Hardy inequalities  (see \cite[Proposition 3]{Miclo1999}) asserts that for any birth-death chain, the optimal constant $C$ in the Poincaré inequality (\ref{BDCPoincaré}) satisfies
\begin{eqnarray}
C & \le & 4\  \inf_{o\in\bI} \left\{B_+(o)\vee B_-(o)\right\},
\end{eqnarray}
where we have introduced the quantities
\begin{eqnarray*}
B_+(o) \ := \ \sup_{x\in(o,\infty)} \sum_{y=o+1}^x\frac{\pi\left(\left[x,\infty\right)\right)}{\pi(y)q_-(y)}  & \textrm{and} &
B_-(o) \ := \ \sup_{x\in(-\infty,o)} \sum_{y=x}^{o-1}\frac{\pi\left(\left(-\infty,x\right]\right)}{\pi(y)q_+(y)}.
\end{eqnarray*}
Here, all intervals  are implicitly restricted to $\bI$ and the supremum over an empty set is interpreted as zero.
Now, for any $x\in\bI$, the detailed balance equation gives
\begin{eqnarray}
\frac{\pi(x+1)}{\pi(x)} & = & \frac{q_+(x)}{q_-(x+1)},
\end{eqnarray}
and our assumption  ensures that the right-hand side is a non-increasing function of $x$. 
Thus, $\pi$ is a log-concave distribution. In particular, for all $o,y,x$ in $\bI$ with $o<y\le x$, we have
\begin{eqnarray}
{\pi(y)} & \ge & \frac{\pi(x)\pi(o)}{\pi(o+x-y)}.
\end{eqnarray}
Using this estimate and $q_-(y)\ge q_-(o+1)\ge q_\star(o)$ inside the definition of $B_+(o)$,  we obtain
\begin{eqnarray*}
B_+(o) & \le & \sup_{x\in(o,\infty)}\left\{\frac{\pi([o,x))\pi([x,\infty))}{\pi(x)\pi(o)q_\star(o)}\right\} \\ & \le & \frac{1}{\pi^2(o)q_\star(o)},
\end{eqnarray*}
where the second inequality uses the fact that $x\mapsto \frac{\pi([x,\infty))}{\pi(x)}$ is non-increasing, by log-concavity of $\pi$.
By symmetry, the same bound applies to $B_-(o)$, and the proof is complete.
\end{proof}
As is clear from its proof, the above lemma is inherently restricted to situations where $\pi$ is log-concave. Fortunately for us, this is the case of Poisson distributions. There are of course many natural choices of monotone birth-death rates that are reversible with respect to a Poisson distribution, and each leads to a different weighted Poincaré inequality. The particular weight appearing in the next proposition will turn out to be the right one for establishing our target IDI in the next section. Let $\phi\colon (0,\infty)\to(0,\infty)$ be defined by
\begin{eqnarray}
\phi(u) & := & \frac{u-1-\log u}{\log^2u},
\end{eqnarray}
with the right-hand side being interpreted as $1/2$ when $u=1$. This function is easily checked to be continuous, increasing and concave. For $u,v>0$, we then introduce the shorthand
\begin{eqnarray}
\label{def:weight}
w(u,v) & := & u\phi\left(\frac{v}{u}\right).
\end{eqnarray}
\begin{proposition}[Weighted Poincaré inequality for the Poisson distribution]\label{pr:poisson}Let $\pi$ be the Poisson distribution with mean $t\in(0,\infty)$. Then, for any $f\in L^2(\pi)$, we have  
\begin{eqnarray}
\Var_{\pi}\left(f\right) & \le & \psi(t)\sum_{x=0}^\infty w\left(\pi(x),\pi(x+1)\right)\left[f(x+1)-f(x)\right]^2,
\end{eqnarray}
where $\psi(t):=550\,t \,\max\left \{1,\log \frac 1t\right\}$.
\end{proposition}
\begin{proof}
Consider the  birth-death chain with stationary law $\pi$ defined as follows: for each $x\in\dN$, the jump rates from $x$ to $x+1$ and from $x+1$ to $x$ are respectively given by
\begin{eqnarray}
q_+(x) & := & \frac{w\left(\pi(x),\pi(x+1)\right)}{\pi(x)} \ = \ \phi \left(\frac{t}{x+1}\right),\\
 q_-(x+1) & := &  \frac{w\left(\pi(x),\pi(x+1)\right)}{\pi(x+1)}\ = \ \frac{x+1}{t}\phi \left(\frac{t}{x+1}\right).
\end{eqnarray}
The fact that $\phi$ is increasing readily implies that  $x\mapsto q_+(x)$ is decreasing. On the other hand, the fact that $\phi$ is concave with $\phi(0+)=0$ ensures that $u\mapsto \phi(u)/u$ is decreasing, so that  $x\mapsto q_-(x)$ is increasing. Thus, Lemma \ref{pr:BDC} applies, and it only remains to show that
\begin{eqnarray}
\label{goal}
\inf_{x\in\dZ_+}\left\{\frac{4}{\pi^2(x)q_\star(x)}\right\} & \le & 550\,t\max\left\{1,\log\frac{1}{t}\right\}.
\end{eqnarray}
To do so, we choose $x=\lfloor t\rfloor$ in the above infimum. For $t\ge 1$,
\begin{eqnarray}
q_\star(\lfloor t\rfloor) & = & \min\left\{\phi\left(\frac{t}{\lfloor t\rfloor}\right),\frac{\lfloor t\rfloor+1}{t}\phi\left(\frac{t}{\lfloor t\rfloor+1}\right)\right\} \ \ge \ \phi(1) \ = \ \frac 12,\\
\pi(\lfloor t\rfloor) & = & \frac{e^{-t}t^{\lfloor t\rfloor}}{\lfloor t\rfloor!} \ \ge \ \frac{1}{3e\sqrt{t}},
\end{eqnarray}
 thanks to the classical estimate $n!\le 3\sqrt{n}(\frac{n}{e})^n$ for $n\ge 1$. On the other hand, for $t<1$,
 \begin{eqnarray}
\pi^2(0)q_\star(0) & =  & \frac{\phi(t)}{te^{2t}} \ \ge \ \frac{1}{3te^2\log \frac 1t} .
\end{eqnarray}
Together, those two estimates guarantee that (\ref{goal}) holds, since $72e^2\approx 532$.
\end{proof}
By tensorization, the above proposition admits the following multi-dimensional extension.
\begin{proposition}[Tensorized version]\label{pr:tensor}Let $\pi=\pi_1\otimes\cdots\otimes\pi_d\in\cP(\dN^d)$ be the tensor product of $d\ge 1$ Poisson laws with respective means $t_1,\ldots,t_d\in(0,\infty)$. Then, for any $f\in L^2(\pi)$, 
\begin{eqnarray}
\Var_{\pi}\left(f\right) & \le & \sum_{i=1}^d \psi(t_i)\sum_{x\in\dN^d} w(\pi(x),\pi(x+e_i))\left[f(x+e_i)-f(x)\right]^2,
\end{eqnarray}
where $(e_1,\ldots,e_n)$ denotes the canonical basis of $\dZ^d$ and $\psi$  the function in Proposition \ref{pr:poisson}.
\end{proposition}
\begin{proof}
Fix $i\in\{1,\ldots,d\}$. Proposition \ref{pr:poisson} allows us to estimate the variance of $f$ when viewed as a function of the $i$-th coordinate only, the  values of the remaining coordinates being fixed. Specifically, for every $(x_1,\ldots,x_{i-1},x_{i+1},\ldots,x_d)\in \dN^{d-1}$, we have
\begin{eqnarray*}
\Var_{\pi_i}\left[f\left(x_1,\ldots,x_{i-1},\cdot,x_{i+1},\ldots,x_d\right)\right] 
& \le & \psi(t_i) \sum_{x_i=0}^\infty w\left(\pi_i(x_i),\pi_i(x_i+1)\right)\left[f(x+e_i)-f(x)\right]^2\\
& = &  \psi(t_i) \sum_{x_i=0}^\infty \pi_i(x_i)\phi\left(\frac{\pi(x+e_i)}{\pi(x)}\right)\left[f(x+e_i)-f(x)\right]^2,
\end{eqnarray*}
where we have used the shorthand $x=(x_1,\ldots,x_d)$ on the right-hand side.
Multiplying by $\pi_1(x_1)\cdots\pi_{i-1}(x_{i-1})\pi_{i+1}(x_{i+1})\cdots\pi_d(x_d)$ and summing over  $x_1,\ldots,x_{i-1},x_{i+1},\ldots,x_d$ yields
\begin{eqnarray*}
\Var_{\pi,i}(f) & \le & \psi(t_i) \sum_{x\in \dN^d}\pi(x)\phi\left(\frac{\pi(x+e_i)}{\pi(x)}\right)\left[f(x+e_i)-f(x)\right]^2,
\end{eqnarray*}
where $\Var_{\pi,i}(f)$ denotes the expected conditional variance of $f$ given all coordinates except the $i-$th one. At this point, the desired conclusion readily follows from the fact that
\begin{eqnarray*}
\Var_{\pi}(f) & \le & \sum_{i=1}^d\Var_{\pi,i}(f),
\end{eqnarray*}
which is the well-known tensorization property of the variance \cite{EfronStein1981}.
\end{proof}
\subsection{Projection onto abelian groups}
As announced long ago, we can finally project our  Poincaré inequality for Poisson distributions onto the abelian group $G$, in order to obtain a \emph{local} Poincaré inequality for the random walk $X=(X_t)_{t\ge 0}$ with jump distribution $\mu$. The fact that our weight has the form 
$w(u,v)  =   u\phi\left(\frac{v}{u}\right)$,
with $\phi$ increasing and concave turns out to play a crucial role here. 
\begin{proposition}[Projected version]\label{pr:proj}Fix $t>0$ and set $\mu_t:=\mathrm{law}(X_t)$. Then,
\begin{eqnarray}
\Var_{\mu_t}\left(f\right) & \le & \Psi(t)\sum_{g,s\in G} \mu(s)w(\mu_t(g),\mu_t(g+s))\left[f(g+s)-f(g)\right]^2,
\end{eqnarray}
for every function $f\colon G\to \dR$, where $\Psi$ is as in Theorem \ref{th:IDI}.
\end{proposition}
\begin{proof}Write $\supp(\mu)=\{s_1,\ldots,s_d\}$ with $d=|\supp(\mu)|$, and let $\pi\in\cP(\dN^d)$ be the tensor product of $d$ Poisson laws with respective means $t_1,\ldots,t_d$, where $t_i:=t\mu(s_i)$. By construction, $\mu_t$ is then the image of $\pi$ under the map $\rho\colon \dN^d\to G$ defined by
\begin{eqnarray}
\rho(v) & := & v_1s_1+\cdots+v_d s_d. 
\end{eqnarray}
Thus, for any $f\in\dR^G$, we have $\Var_{\mu_t}\left(f\right) =  \Var_{\pi}(f\circ\rho)$ and  Proposition \ref{pr:tensor} gives 
\begin{eqnarray*}
\Var_{\mu_t}\left(f\right) & \le & \sum_{i=1}^d \psi(t_i)\sum_{x\in\dN^d} w\left(\pi(x),\pi(x+e_i)\right)\left[(f\circ \rho)(x+e_i)-(f\circ\rho)(x)\right]^2\\
& = & \sum_{i=1}^d \psi(t_i)\sum_{g\in G}\sum_{x\in\rho^{-1}(g)} w\left(\pi(x),\pi(x+e_i)\right)\left[f(g+s_i)-f(g)\right]^2\\
& \le & \sum_{i=1}^d \psi(t_i)\sum_{g\in G} w\left(\sum_{x\in\rho^{-1}(g)}\pi(x),\sum_{x\in\rho^{-1}(g)}\pi(x+e_i)\right)\left[f(g+s_i)-f(g)\right]^2,
\end{eqnarray*}
where we have crucially used the super-additivity property
\begin{eqnarray}
\label{superadd}
\sum_{i\in I}w\left(u_i,v_i\right) \ = \ \sum_{i\in I}u_i\phi\left(\frac{v_i}{u_i}\right)& \le & \sum_{i\in I}u_i\phi\left(\frac{\sum_{i\in I}v_i}{\sum_{i\in I}u_i}\right) \ = \  w\left(\sum_{i\in I}u_i,\sum_{i\in I}v_i\right),
\end{eqnarray}
which holds for any countable collection  $(u_i,v_i)_{i\in I}$ of positive numbers, by  concavity of $\phi$. Finally, since $\mu_t$ is the image of $\pi$ under $\rho$, we have for all $g\in G$,
\begin{eqnarray}
\sum_{x\in\rho^{-1}(g)}\pi(x) \ = \ \mu_t(g) & \textrm{and} &  \sum_{x\in\rho^{-1}(g)}\pi(x+e_i) \ \le \ \mu_t(g+s_i).
\end{eqnarray}
Since $w(\cdot,\cdot)$ is non-decreasing in its second argument, this implies 
\begin{eqnarray}
w\left(\sum_{x\in\rho^{-1}(g)}\pi(x),\sum_{x\in\rho^{-1}(g)}\pi(x+e_i)\right) & \le &  w\left(\mu_t(g),\mu_t(g+s_i)\right).
\end{eqnarray}
Inserting this into the previous computation concludes the proof, since $\psi(t_i) \le \Psi(t)\mu(s_i)$.
\end{proof}
We are now in a position to establish the universal IDI promised in Theorem \ref{th:IDI}.
\begin{proof}[Proof of Theorem \ref{th:IDI}]Fix $t>0$ and set $\mu_t:=\mathrm{law}(X_t)$. Applying Proposition \ref{pr:proj} to the function $f\colon g\mapsto \log \mu_t(g)$, and using the particular expression of $w$, we obtain
\begin{eqnarray*}
\Varent(X_t) & \le &  \Psi(t)\sum_{s\in G} \mu(s)\sum_{g\in G} w(\mu_t(g),\mu_t(g+s))\left[\log \frac{\mu_t(g+s)}{\mu_t(g)}\right]^2\\
& = &   \Psi(t)\sum_{s\in G} \mu(s)\sum_{g\in G} \left[\mu_t(g+s)-\mu_t(g)-\mu_t(g)\log \frac{\mu_t(g+s)}{\mu_t(g)}\right]\\
& = &  \Psi(t)\frac{\dd}{\dd t}\Ent(X_t),
\end{eqnarray*}
where the last line classically follows from the Fokker$-$Plank equation $\frac{\dd}{\dd t}\mu_t=\cL^\top\mu_t$.
\end{proof}
 We have already explained how this leads to our main estimate (Theorem \ref{th:wmix}), and the latter readily implies our general sufficient criterion for cutoff (Corollary \ref{co:cutoff}). Finally, the fact that the latter reduces to the product condition in the logarithmically sparse regime (Corollary \ref{co:sparse}) is a direct consequence of the observation that
\begin{eqnarray}
\left|\supp(\mu)\right| & \le & \frac{1}{\mu_\star},
\end{eqnarray}
together with the following simple  estimate on the mixing time of abelian random walks.
\begin{lemma}[Slow mixing for sparse abelian random walks]\label{lm:LB}For every $a>0$, there is a  constant $\eta=\eta(a)>0$, independent of the pair $(G,\mu)$, such that
\begin{eqnarray}
\left|\supp(\mu)\right| \ \le \ a\log |G| & \Longrightarrow & \tmix \ \ge \ \eta\log |G|.
\end{eqnarray}
\end{lemma}
\begin{proof}
We can assume that $|G|\ge 2$, otherwise the claim is trivial. Write $\supp(\mu)=\{s_1,\ldots,s_d\}$ with $d=|\supp(\mu)|$, and fix  $t\ge 0$. By Markov's inequality, there is a chance at least $3/4$  that  no more than $r=\lfloor 4t\rfloor$ jumps are made during the time interval $[0,t]$. Moreover, on that event, the random variable $X_t$ is confined to the set 
\begin{eqnarray}
A & := &  \left\{x_1 s_1+\cdots+x_ds_d\colon (x_1,\ldots,x_d)\in\dN^d, x_1+\cdots+x_d\le r \right\}.
\end{eqnarray}
This combinatorial description implies that 
\begin{eqnarray}
|A|  \ \le  \ { r+d\choose  r} \ \le \ e^{h(r,d)}, &\textrm{ where }& h(r,d)\ :=\ r+r\log\left(1+\frac{d}{r}\right).
\end{eqnarray}
Note that the function $h\colon(0,\infty)^2\to(0,\infty)$ thus defined is increasing. In particular, if
\begin{eqnarray}
d \ \le \ a\log |G|, & \textrm{ and }& t \ \le \ \eta \log|G|
\end{eqnarray}
for some $a,\eta\in(0,\infty)$, then recalling that $r=\lfloor 4t\rfloor$, we deduce that
\begin{eqnarray}
|A| & \le & e^{h(4\eta\log|G|,a\log |G|)} \ = \ |G|^{h(4\eta,a)}.
\end{eqnarray}
Putting things together, we conclude that 
\begin{eqnarray}
\label{LB}
\dd(t) & \ge & \PP(X_t\in A)-\cU(A) \ \ge \ \frac{3}{4}-\left(\frac{1}{|G|}\right)^{1-h(4\eta,a)}.
\end{eqnarray}
Finally, note that once $a>0$ is fixed, the exponent $1-h(4\eta,a)$ can be made arbitrarily close to $1$ by choosing $\eta$ small enough. In particular, it can be chosen so that the right-hand side of (\ref{LB}) is larger than $1/(2e)$,  independently of the pair $(G,\mu)$.
\end{proof}
\section{Proof of sharpness}
Our proof of Theorem \ref{th:sharp} builds upon the  classical observation that  rank-one perturbations can destroy cutoff while preserving the product condition. More precisely, consider a finite group $G$, a jump distribution $\mu\in\cP(G)$ and a number $\theta\in(0,1)$, all depending on an implicit parameter $n\ge 1$ in such a way that the random walk with jump distribution $\mu$ satisfies $\gap\,\tmix\gg 1$ and $\theta\,\tmix=\Theta(1)$. Then the random walk with jump distribution 
\begin{eqnarray}
\label{rank-one}
\nu & := & (1-\theta)\mu+\theta \cU 
\end{eqnarray}
  still satisfies  the product condition, but fails to exhibit cutoff  (see, e.g., \cite[Section 1.5]{salez2025lecture}).  Unfortunately, such  examples are not \emph{reasonable}, in the sense that the  jump distribution (\ref{rank-one}) has full support. To obtain sparser examples, a natural idea is to replace the uniform law  by a \emph{sparse approximation} of it, i.e. a measure $\sigma\in\cP(G)$ whose support is much smaller, but which still mixes extremely fast. Here, in addition, we would like the resulting  measure 
\begin{eqnarray}
\nu & := & (1-\theta)\mu+\theta\sigma 
  \end{eqnarray}  
to be uniform on its support. This  forces $\mu$ and $\sigma$ to be uniform on disjoint sets, with sizes proportional to $1-\theta$ and $\theta$ respectively. Recalling that our initial jump distribution $\mu$ should satisfy the product condition and have mixing time of order $1/\theta$, we obtain a rather constrained set of specifications. The latter turns out to be feasible, as we will see. 
 To keep our analysis as simple as possible,  we implement this strategy on cyclic groups. 
\subsection{Fourier analysis on cyclic groups}The \emph{Fourier transform} of a probability measure $\mu$ on the cyclic group $G=\dZ_n$ is given by
\begin{eqnarray}
\label{def:FT}
{\wnu}(z)   & := & \sum_{y\in G}\mu(y)\exp\left(\frac{2i\pi zy}{n}\right),
\end{eqnarray}
for $z\in G$. We define the \emph{spectral radius} $\varrho=\varrho(\mu)$ as the $[0,1]$-valued number
\begin{eqnarray}
\varrho & := & \max_{z\in G\setminus\{0\}}|\wnu(z)|.
\end{eqnarray}
This parameter  controls  the mixing properties of the random walk, as summarized next. 
\begin{lemma}[Spectral bounds]\label{lm:UB}The random walk with jump distribution $\mu$ satisfies
\begin{eqnarray}
\gap \ \ge \ 1-\varrho &\textrm{ and } & 
\dd(t) \ \le \ {\sqrt{n} \varrho^{r+1}}+ \sum_{k=0}^{r}p_t(k),
\end{eqnarray}
for all $t\ge 0$ and all $r\in\dN$, where $p_t\in\cP(\dN)$ denotes the Poisson law with mean $t$.
\end{lemma}
\begin{proof}
The transition matrix $T(x,y):=\mu(y-x)$ describing the effect of a jump according to $\mu$ can readily be diagonalized in the standard Fourier basis of $L^2(G)$, and its eigenvalues are exactly the complex numbers $\wnu(z)$, $z\in G$. Since the generator (\ref{def:L}) is   $\cL=T-{\rm Id}$, the first claim follows. Moreover, the Cauchy-Schwarz inequality yields 
\begin{eqnarray}
2\dtv(\mu,\cU) & = & \sum_{x\in G}\left|\mu(x)-\frac{1}{n}\right|
\\ & \le & \sqrt{n\sum_{x\in G}\left|\mu(x)-\frac{1}{n}\right|^2} \\
 & = &   \sqrt{\sum_{z\in G\setminus\{0\}}\left|\widehat{\mu}(z)\right|^2}  \\
 & \le & \sqrt{n}\varrho,
\end{eqnarray}
where the third line is   Parseval's identity. Replacing $\mu$ with its $k-$fold convolution $\mu^{\star k}$ and using the convolution-product property 
$
\widehat{\mu_1\star\mu_2}  =  \widehat{\mu_1}\,\widehat{\mu_2}
$ of the Fourier transform, we obtain
\begin{eqnarray}
\dtv(\mu^{\star k},\cU) & \le &   {\sqrt{n}\varrho^k},
\end{eqnarray}
for each $k\in\dN$. Finally, note that by definition, our random walk $X=(X_t)_{t\ge 0}$ satisfies
\begin{eqnarray}
\mathrm{law}(X_t) & = & \sum_{k=0}^\infty p_t(k)\mu^{\star k}.
\end{eqnarray}
Since the total-variation distance is convex, we conclude that 
\begin{eqnarray}
\label{prelim:UB}
\dd(t) & \le & \sum_{k=0}^\infty  p_t(k)\dtv(\mu^{\star k},\cU)\\
& \le &  {\sqrt{n}\varrho^{r+1}}+ \sum_{k=0}^{r}p_t(k), 
\end{eqnarray}
for any $t\ge 0$ and any $r\in\dN$, as desired.
\end{proof}
\begin{lemma}[Random measures mix quickly]\label{lm:random}Consider the (symmetric) random measure 
\begin{eqnarray}
\label{def:random}
\mu & := & \frac{1}{d}\sum_{k=1}^d\frac{\delta_{Y_k}+\delta_{-Y_k}}{2},
\end{eqnarray}
where $d\ge 1$ is an integer and $Y_1,\ldots,Y_d$ are independent and uniform on $G=\dZ_n$. 
Then, 
\begin{eqnarray}
\PP\left(\varrho(\mu)\ge 2\sqrt{\frac{\log n}{d}}\right) & \le & \frac{1}{n}.
\end{eqnarray}
\end{lemma}
\begin{proof}
Fix $z\in G\setminus\{0\}$. When applied to the random measure (\ref{def:random}), the formula (\ref{def:FT}) gives
\begin{eqnarray}
\wnu(z) & =  &  \frac{1}{d}\sum_{k=1}^d \cos\left(\frac{2\pi z Y_k}{n}\right).
\end{eqnarray}
This is the empirical mean of $d$ independent, centered, $[-1,1]$-valued variables, so 
\begin{eqnarray}
\PP\left(\left|\wnu(z)\right|\ge a\right) & \le & 2e^{-\frac{da^2}{2}},
\end{eqnarray}
for any $a\ge 0$, by Hoeffding's inequality. Choosing $a=2\sqrt{(\log n)/d}$ sets the right-hand  side to $2/n^2$. Since $\wnu(-z)=\wnu(z)$, taking a union bound over half of the $z\in G\setminus\{0\}$ suffices.
\end{proof}

\subsection{Achieving polynomial degrees}
The above lemmas already allow us to prove the following dense version of Theorem \ref{th:sharp}. 
\begin{proposition}[Polynomially sparse examples]\label{pr:poly}Fix $\alpha\in(0,1/2)$. Then, for each $n\ge 1$, there is a symmetric set $S_n\subseteq\dZ_n$ such that simple random walk on $\mathrm{Cay}(\dZ_n,S_n)$ satisfies
\begin{enumerate}[(i)]
\item $|S_n|\sim n^\alpha$;
\item $\gap=\Omega(1)$;
\item $\tmix=\Theta(\log n)$;
\item $\displaystyle{\limsup_{n\to\infty}\dd(u\log n)<1}$ for every $u>0$.
\end{enumerate}
In particular, the product condition is satisfied but there is no cutoff.
\end{proposition}
\begin{proof}We omit the index $n$ to lighten notation. 
Set $d:=\lfloor 10\log n\rfloor$ and $m:=\lfloor n^\alpha/(40\log n)\rfloor$. Let $Y_1,\ldots,Y_d,Z_1,\ldots,Z_m$ be independent and uniformly distributed on $\dZ_n$ and consider the random symmetric subset of $\dZ_n$ defined as follows:
\begin{eqnarray}
\label{def:S}
S & := & \pm\left(\{Z_1,\ldots, Z_{m}\}\cup \bigcup_{i=1}^d [Y_i-m,Y_i+m]\right).
\end{eqnarray}
We claim that $S$ has the desired properties with high probability as $n\to\infty$. First observe that $(dm)^2\ll n$, so that with high probability, no repetition occurs in the list of elements appearing on the right-hand side of (\ref{def:S}). On that event, we have
\begin{eqnarray}
|S| & = & 2m+2d(2m+1) \ \sim \ n^\alpha,
\end{eqnarray}
showing (i). Moreover, on the same event, the uniform law on $S$ coincides with 
\begin{eqnarray}
\label{def:nu}
\nu & := & (1-\theta)(\xi\star\mu)+\theta \sigma,
\end{eqnarray}
where $\theta=\frac{2m}{|S|}\sim \frac{1}{20\log n}$ and  $\xi$ is the uniform measure on $[-m,m]$, while
\begin{eqnarray}
\mu \ := \ \frac{1}{2d}\sum_{j=1}^d(\delta_{Y_j}+\delta_{-Y_j})  & \textrm{ and } &  \sigma \ := \ \frac{1}{2m}\sum_{j=1}^m(\delta_{Z_j}+\delta_{-Z_j}).
\end{eqnarray}
By Lemma \ref{lm:random}, we have with high probability,
\begin{eqnarray}
\label{ctx:markov}
\varrho(\mu)  \ \le \ \frac 23 & \textrm{ and } & \varrho(\sigma)  \  \le \  \frac{13\log n}{n^{\frac\alpha 2}}.
\end{eqnarray}
By  linearity and the convolution-product property of the Fourier transform, this implies $\varrho(\nu)\le 2/3$ which, in view of Lemma \ref{lm:UB}, guarantees (ii), as well as half of (iii), namely $\tmix=O(\log n)$. To prove the other half, fix $u>0$  and set $t:=\lceil u\log n\rceil$. During the time interval $[0,t]$, a random walk with jump distribution $\nu$ has chance 
\begin{eqnarray}
q & := & e^{-\theta t}\sum_{k=0}^{t} p_{(1-\theta)t}(k)\ \sim \ \frac{1}{2}e^{-\frac{u}{20}}
\end{eqnarray}
to use no jump in $\supp(\sigma)$ and no more than $t$ jumps in $\supp(\xi\star \mu)$, in which case it lies in
\begin{eqnarray}
A & := & [-tm,tm]+\bigcup_{i=0}^t\supp(\mu^{\star i}).
\end{eqnarray}
Using the function $h$ in the proof of Lemma \ref{lm:LB}, we have 
\begin{eqnarray}
|A| & \le & (2tm+1){2d+t\choose t} \ = \ O\left(n^{\frac 12+h(u,20)}\right).
\end{eqnarray}
Choosing $u$ so small that $h(u,20)<1/2$ ensures that $\cU(A)\ll 1$, and we conclude that 
\begin{eqnarray}
\liminf_{n\to\infty}\dd(u\log n) & \ge & \frac{1}{2}e^{-\frac{u}{20}} \ \ge \ \frac 14,
\end{eqnarray}
establishing $\tmix=\Omega(\log n)$. Finally, since the distance to equilibrium decreases under convolution, $\dd(t)$ is at most the distance to equilibrium at time $\theta t$ for a random walk with jump distribution $\sigma$, and applying Lemma \ref{lm:UB} to the latter yields
\begin{eqnarray*}
\dd(t) 
& \le & {\sqrt{n}}[\varrho(\sigma)]^{r+1}+  \sum_{j=0}^{r}p_{\theta t}(j),
\end{eqnarray*}
 for all $t\ge 0$ and $r\in\dN$. In view of (\ref{ctx:markov}), choosing $r=\lfloor 1/\alpha\rfloor$ is enough to guarantee that the first term vanishes, and since $\theta\sim 1/(20\log n)$, we conclude that for every $u>0$,
 \begin{eqnarray*}
\limsup_{n\to\infty}\dd(u\log n) & \le &  \sum_{j=0}^{r}p_{u/20}(j). 
\end{eqnarray*}
The right-hand side is of course less than $1$, and (iv) is proved. 
\end{proof}
\subsection{Achieving near-logarithmic degrees}In this final section, we improve the polynomial sparsity of the above construction to the near-logarithmic one promised in Theorem \ref{th:sharp}. Rather than decreasing the degree of the graph, our idea is to artificially increase its volume by embedding the walk into a much larger group. Specifically, we add an independent coordinate  whose dimension is exponentially large, but which mixes so fast that it does not affect the overall distance to equilibrium. 
\begin{proof}[Proof of Theorem \ref{th:sharp}]We use the notation $\gap(\Gamma), \tmix(\Gamma)$ and $\dd_\Gamma(\cdot)$ for the spectral gap, mixing time and distance to equilibrium of a simple random walk on the abelian Cayley graph  $\Gamma$.  Let $(a_n)_{n\ge 1}$ satisfy $1\ll \log a_n\ll\log n$. We consider two ($n$-dependent) symmetric abelian Cayley graphs $\Gamma_1=\mathrm{Cay}(G_1,S_1)$ and $\Gamma_2=\mathrm{Cay}(G_2,S_2)$. The first one is such that
\begin{eqnarray}
|G_1|=n,\quad |S_1|\sim  a_n\log n,\quad \gap(\Gamma_1)  = \Theta(1),\quad \tmix(\Gamma_1) \ll \log n,
\end{eqnarray}
and its existence is guaranteed, for example, by the random construction in Lemma \ref{lm:random} (take $G_1=\dZ_n$ and $S_1=\{\pm Y_1,\ldots,\pm Y_{d}\}$ with  $d=\lfloor (a_n \log n)/2\rfloor$ and $Y_1,\ldots,Y_d$ independent and uniform). For the second one, we simply invoke Proposition \ref{pr:poly} with $\alpha=1/3$ to obtain, after extraction, a graph satisfying
\begin{eqnarray}
|G_2|\sim a_n^3,\quad |S_2|\sim  a_n,\quad \gap(\Gamma_2)  =  \Theta(1),\quad \tmix(\Gamma_2) = \Theta(\log a_n),
\end{eqnarray}
as well as the following  property, which crucially precludes the occurrence of a cutoff:
\begin{eqnarray}
\label{eq:nocutoff}
\forall u>0,\qquad \limsup_{n\to\infty}\dd_{\Gamma_2}(u\log a_n) & < &  1.
\end{eqnarray}
Now, we claim that the properties announced in the theorem are satisfied by the \emph{Cartesian product}  of  $\Gamma_1$ and $\Gamma_2$, i.e. the Cayley graph $\Gamma=\mathrm{Cay}(G,S)$ defined by
\begin{eqnarray}
G \ := \ G_1\times G_2, & \textrm{ and } & S\ :=\ (S_1\times\{0\})\cup (\{0\}\times S_2).
\end{eqnarray}
We first note that this graph has the desired sparsity $|S|\sim a_n\log|G|$. Indeed,
\begin{eqnarray}
\log |G| \ = \ \log |G_1|+\log |G_2| \ \sim \ \log n & \textrm{ and } & 
|S| = |S_1|+|S_2| \ \sim \ a_n\log n.
\end{eqnarray}
Also, observe that the simple random walk on $\Gamma$ can be described as $X(t) =  \left(X_1(\theta_1 t),X_2(\theta_2t)\right)$, 
where $X_1,X_2$ are independent simple random walks on $\Gamma_1,\Gamma_2$ respectively, and where
\begin{eqnarray}
\theta_1 \ = \ \frac{|S_1|}{|S|} \ \sim 1 & \textrm{ and } & \theta_2 \ = \ \frac{|S_2|}{|S|} \ \sim \frac{1}{\log n}.
\end{eqnarray}
This decomposition into independent coordinates guarantees that 
\begin{eqnarray}
\max\left\{\dd_{\Gamma_1}(\theta_1 t),\dd_{\Gamma_2}(\theta_2 t)\right\} \ \le \ \dd_\Gamma(t)\ \le \ \dd_{\Gamma_1}(\theta_1 t)+\dd_{\Gamma_2}(\theta_2 t),
\end{eqnarray}
for all $t\ge 0$. In particular, it follows that
\begin{eqnarray}
\tmix(\Gamma) & = & \Theta\left(\frac{\tmix(\Gamma_1)}{\theta_1}+\frac{\tmix(\Gamma_2)}{\theta_2}\right) \ = \ \Theta(\log n\log a_n),
\end{eqnarray}
while (\ref{eq:nocutoff}) shows that for any fixed $u>0$,
\begin{eqnarray}
\limsup_{n\to\infty}\dd_\Gamma\left(u\tmix(\Gamma)\right) & < & 1.
\end{eqnarray}
Thus, the simple random walk on $\Gamma$ does not exhibit cutoff. Yet, we have
\begin{eqnarray}
\gap(\Gamma) \ = \ \min\left\{\theta_1 \gap(\Gamma_1),\theta_2 \gap(\Gamma_2)\right\} \ = \ \Theta\left(\frac 1{\log n}\right),
\end{eqnarray}
so that the product condition is satisfied. Finally, the absence of cutoff implies that  for small enough $\varepsilon>0$,  $\wmix(\varepsilon)=\Theta(\tmix)=\Theta(\log n\log a_n)$, while a direct evaluation shows that the right-hand side of our main window estimate (Theorem \ref{th:wmix}) is also $\Theta(\log n\log a_n)$. 
\end{proof}
\section*{Acknowledgment and AI-use statement}
This work was partially supported by the ERC consolidator grant 101123174. AI was only used for typos and bibliography. The writing, ideas and proofs are entirely our own. 
\bibliographystyle{plain}
\bibliography{draft}

\end{document}